**Intuitive, didactically useful, and historically possible proofs for the two Egyptian pyramid volume "formulas" (1850 BCE): Thoughts on the border between history and didactics of mathematics**

Reinhard Siegmund-Schultze (University of Agder, Kristiansand, Norway)[1]

**Abstract:**

Egyptologists and historians of mathematics around 1930 did an admirable job in showing that problem 14 of the newly discovered Moscow Papyrus from around 1850 BCE amounts to a general and exact calculation of the volume of a truncated pyramid (frustum). They were less successful in giving tentative explanations how the Egyptians may have found the solution and what convinced them of its correctness. In particular, those historians never looked into the possibility of dissecting three identical copies of the frustum and recomposing them into three boxes of differing sizes whose volumes can be easily calculated. This is surprising because the "formula" at which the historians arrived seems to suggest this procedure. About 2000 years after the Egyptians, the Chinese scholar Liu Hui did exactly this for the almost identical problem from the "Nine Chapters". If those "historians of mathematics around 1930" had known Liu Hui's algorithm they could have easily drawn tentative conclusions also for the Egyptian case. The present paper suggests that it was their knowledge of rigorous Euclidean geometry and of relatively "modern" algebra which distorted the judgment of those historians, something which could not have been the case for Liu Hui. Their failure seems to have discouraged later Egyptologists to look for explanations or to even point to Liu Hui when his work had become known. Thus a chance was missed to use the great intuitive and pedagogic potential of a remarkable piece of Egyptian mathematics. This paper is not a contribution to the historiography of Egyptian mathematics for which the author is no specialist but argues primarily on a methodological level using secondary historical sources. The paper is partly inspired by a more recent publication of Paul Shutler (2009) who suggests an intuitive and historically possible proof also for the special, and mathematically crucial case of the "second formula", which is related to the full pyramid.

1 A much shortened version of this paper is about to appear in "British Journal for the History of Mathematics" (2022) under the title: "Another look at the two Egyptian pyramid volume "formulas" of 1850 BCE".

## 1. Introduction and aims

Egyptologists and historians of mathematics around 1930 such as Struve, Gunn, Peet and Neugebauer gave a convincing interpretation of problem 14 of the Moscow Papyrus, found a few decades before, as providing a "formula" for the volume $V_T$ of the truncated pyramid with quadratic base (henceforth "frustum," which is a general notion for a truncated solid). They expressed this algebraically by

$$V_T = \frac{h}{3}\,(a^2 + ab + b^2) \qquad \text{(henceforth } F_T\text{)},$$

with $a$ and $b$ being the sides of the base and the top squares respectively and $h$ being the height of the frustum.

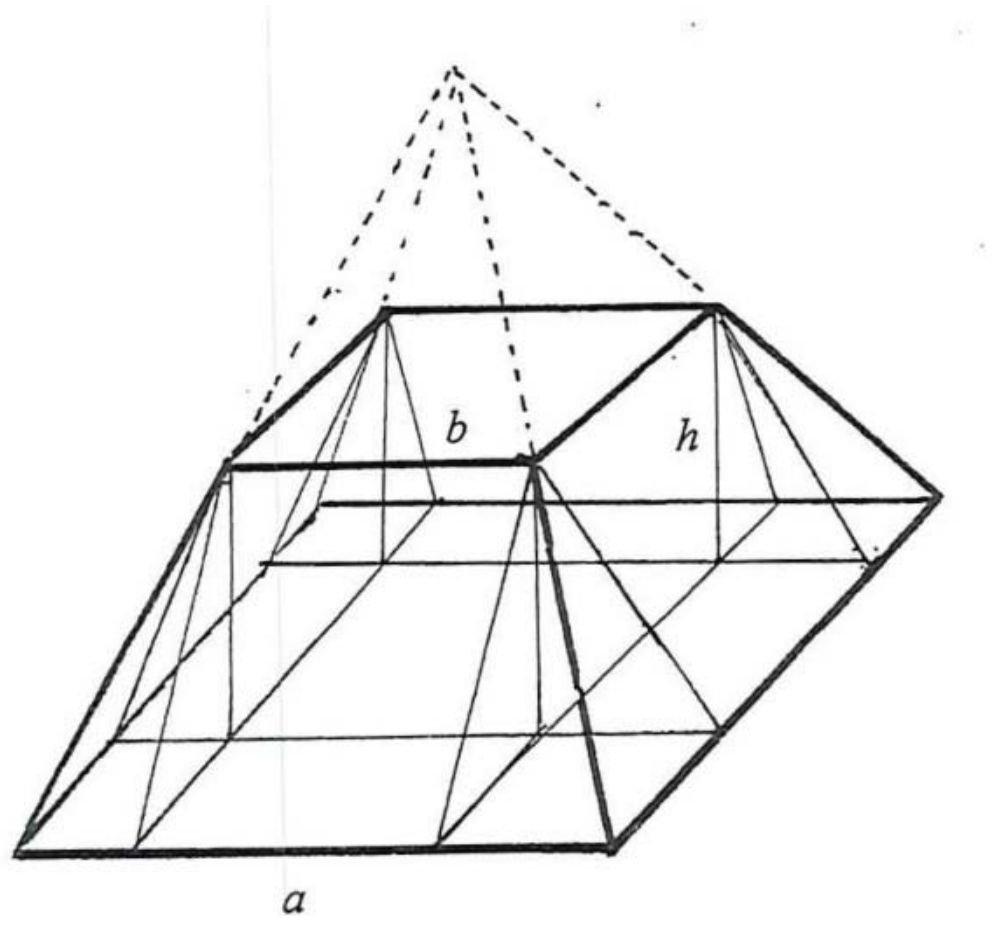

**Figure 1:** A symmetrical truncated pyramid (frustum)

The same historians argued for a variety of empirical, intuitive and somewhat more formal methods which might have led the Egyptians to their formula. Among these were weighing models of the frustum, dissecting it into a cuboid + four prisms + four corner pyramids (as clearly visible from figure 1), understanding the formula as an average of either three volumes or of three areas (to be multiplied afterwards by a common factor $h$), acknowledging certain simple arithmetical operations such as taking out a common factor $h/3$. None of these works to my knowledge considered the formula as expressing the rearrangement of three identical copies of the frustum into three boxes assuming the demolition of the 12 corner pyramids and the reuse of their volume (material) for building up four corner cuboids in the biggest of the resulting three boxes (box (d) in figure 2 below). A few decades later, in the 1970s, scholars of Chinese mathematics showed that this method was basically the one which Liu Hui had used in the third century CE in order to confirm an analogous "formula" from the classical *Nine Chapters* of Chinese mathematics. The method is expressed in the following figure from a paper of Shutler (2009).

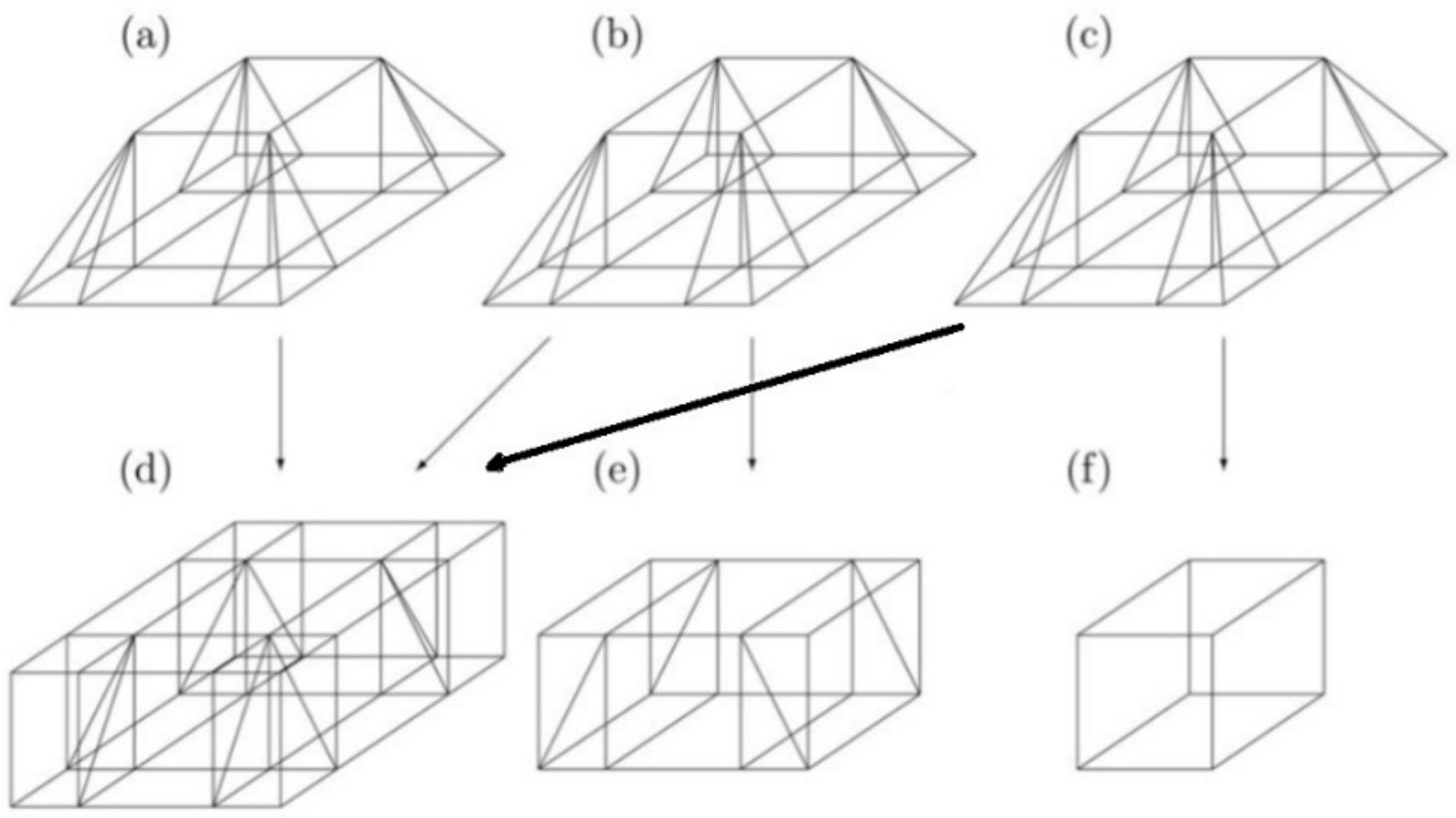

**Figure 2:** from Shutler (2009, 349). I have added the bold arrow from (c) to (d) to Shutler's picture, replacing an erroneous arrow in Shutler which goes from (c) to (e). An explanation follows in section 6 below.

Interestingly, Shutler did not mention Liu Hui or any literature referring to him. My conjecture[2] is that Liu Hui's procedure was in the back of his mind from previous reading, while it was not mentioned in any literature on Egyptian mathematics available to him.

The main stimulus of my following remarks comes from this paper of Shutler's, who – in addition to his "rediscovery" of Liu Hui – gave in it an intuitive and convincing derivation of the formula for the volume of the full pyramid which enables the "demolition"-argument mentioned above. Thus this paper also aims at making Shutler's valuable publication a bit more visible and known.

For my part I will add some methodological reflections concerning the historiography of mathematics and its use for pedagogical purposes.

Deriving formula $F_T$ from three copies of the frustum seemed so natural to me and so much suggested by the formula itself that I had "found" the same method long before having read Shutler, probably though also with Liu Hui in the back of my mind.

[2] I have not succeeded in contacting Shutler per email or otherwise.

Given the broad variety of methods ascribed to the Egyptians by the "historians around 1930"[3] I was wondering why they – not knowing at the time about Liu Hui – had not come up with the same idea as him. Was Liu Hui's method not so trivial and obvious after all as it seemed to me now? Did older (around 1930) or modern historians of Egyptian mathematics find an approach in the style of Liu Hui's somehow incompatible with Egyptian thinking? Were those historians around 1930 influenced in their approach by modern methodology and mathematical knowledge which had not been available to Liu Hui himself in the third century CE? In any case it seemed to me that there was an opportunity – so far barely used – to compare the two very different periods and cultural domains, the Egyptian and the Chinese under one specific aspect and to the benefit of modern students thus using the didactical and intuitive appeal which the construction of pyramids bears for mathematical teaching.

This paper is not a contribution to the historiography of Egyptian mathematics. It is written from the standpoint of a historian of mathematics mostly dealing with more modern periods and with no specialist knowledge of Egyptian culture or language. My didactical reflections are not based on theoretical works or training in this field but on the naïve assumption that ancient Egyptian mathematicians and modern school children with only rudimentary mathematical education have "something in common" when looking at the impressive Egyptian pyramids as mathematical objects.

## 2. The "two pyramid formulas" in Egyptian mathematics

Egyptologists agree that for the Egyptians the construction of their wonderful symmetrical pyramids was one of the most important engineering tasks, that this activity went on through centuries, caused great costs, involved large parts of the population, and required exact planning and calculation of the necessary material and the time to be spent. It must therefore have been known as a "fact" (based on practical work) to those responsible that the volume $V_P$ of the pyramid contained (about) one third of the volume of the (square) cuboid (henceforth "cuboid") with the same base area $a^2$ and the same height $h$. Mathematically written, this means the following formula, which I call henceforth the "implicit Egyptian pyramid volume formula" or "$\frac{1}{3}$ – formula" $F_P$:

$$V_P = \frac{h}{3}\, a^2 \qquad (F_P)$$

[3] I will henceforth use this collective expression to mean the first interpreters of the Moscow Papyrus, including Egyptologists and mathematicians. I have particularly, but not exclusively, in mind (Gunn/Peet 1929), (Struve 1930), and (Neugebauer 1933, 1934).

Similar to physics, the experiences and observations during the construction of the pyramids could be improved to an even higher level of accuracy by experiments with models. The factor of 3 or $\frac{1}{3}$ for the calculation of the pyramids could probably be confirmed by weighing out models of the pyramids, which for example could be made from dried Nile mud, as Gunn and Peet assumed in 1929.[4] Moreover, for special cases such as pyramids inscribed in cubes and in half cubes[5] elementary dissections exist which provide exact mathematical confirmation for the formula in these cases. But we do not know whether the Egyptians knew these dissections.[6] Intuitive arguments[7] (stretching of the cube with unchanged base area in vertical direction) would then allow some semi-mathematical confirmation of the formula for the general case. However, in my opinion this "stretching argument" would not have been intuitive enough for the old Egyptians. In section 8 I will report – based on Shutler (2009) – on a more plausible intuitive proof of $F_P$ which connects directly to the Egyptian example in the Moscow Papyrus.

However, this elementary mathematical formula $F_P$ with great importance for applications, which mathematicians of later centuries extended to asymmetrical pyramids, to pyramids with arbitrary base areas, including circular cones, does not *explicitly* appear in the few surviving Egyptian documents (mostly papyri). This could be due to the fact that there was no actual need to include the formula – because of its simplicity and ease of memorization – among the rules or algorithms apparently written down for the "engineers" or workers. By comparison, in order to calculate circular areas, it was necessary for the

---

[4] In (Struve 1930, 174/75) measuring of the volume of models by immersing them into water is discussed as another possibility.
[5] The "Juel pyramid" with height $a/2$, where the slope of the triangular faces is $45^0$, is one sixth of an $a^3$-cube. The apexes of six identical Juel pyramids meet in the center of the cube. See (Gillings 1964). It is apparently named after the Danish geometer Christian S. Juel (1855–1935). Jesper Lützen kindly alerted me to (Juel 1903) which mentions and partly discusses this type of pyramid.
[6] The dissection of the cube into six Juel pyramids seems more intuitive than the dissection of the cube into three asymmetrical pyramids which was known to the Chinese (Liu Hui). See below figure 4 and sections 4 and 6.
[7] The notion of "intuition" is of course relative and historically dependent. In mathematics there exists intuition on many levels, depending on the educational background, certainly not restricted to geometrical intuition. The notion should not be misunderstood as being opposed to rigor either. I will argue below that Liu Hui's proof of the volume formula for the truncated pyramid was intuitive and rigorous at the same time.

Egyptians to specify and to fix constants which in our interpretation correspond to the rather good value $\frac{256}{81} = 3 + \frac{13}{81} = 3.1604 \ldots$ for π (Imhausen 2016, 118ff.), and apparently therefore this more complicated rule had to be written on some papyri. But, of course, we do not know whether there were other papyri which did not survive or can still be recovered.

In the late 19th century archeologists found the so-called Moscow Papyrus from about the year 1850 BCE where problem 14 contains (after proper interpretation by historians) the formula $F_T$ of the truncated pyramid (frustum), as introduced above.

The latter formula is commonly considered as a highlight of Egyptian mathematics, although it is acknowledged that the available documents show no trace of "proof" or any other weaker mathematical argument. From $F_T$ follows $F_P$ for $b = 0$. Thus we have it finally – if implicitly – on paper (papyrus), a fact of which there was no doubt, namely that the Egyptians knew the $\frac{1}{3}$ – formula $F_P$ for the full pyramid albeit possibly without proof. We are thus allowed (as above in the title) to speak about *two* Egyptian formulas connected to pyramid volumes, one *explicit* $F_T$ and one *implicit* $F_P$.

With today's mathematics both formulas can be easily derived by elementary geometrical, algebraic and analytical means. The historical problem is that only elementary and intuitive geometrical methods can be safely assumed to have been in the possession of the Egyptians, not, however, algebraic or analytical ones. Although one can argue that – mathematically – the core of the problem lies in formula $F_P$ , it is – from a historical (!) standpoint – no triviality to derive $F_T$ from $F_P$. The assumption that $F_T$ could be easily found by subtracting a top pyramid from a bigger one is erroneous because the height of the bigger pyramid has to be determined first (for instance by similarity considerations) which is not viable for a formula needed in applications. Moreover, to derive $F_T$ theoretically from two different $F_P$ is a simple, but non-trivial algebraic assignment even for us today.[8] However, it seems to me that some historians of Egyptian mathematics around 1930 did not sufficiently stress the historical difference between the two formulas, which are in a sense "mathematically equivalent". In particular, W. Struve's mathematical commentary on the Egyptian formula $F_T$ uses all kinds of knowledge about proportions in similar triangles which cannot be assumed to have been in reach of the Egyptians (Struve 1930, 174-176).[9]

[8] This is exactly what Struve does in his mathematical commentary (Struve 1930, 174-176).

[9] See the remarks in (Imhausen 2016, 5) about "anachronistic modern terminology" which are related to some of those historians around 1930.

Given the lack of more extant documents we cannot be sure whether the Egyptians had – beyond experience and possible experimenting – any further mathematical explanations – even intuitive ones – for their formulas at all. Below I will mostly rely on the opinions of specialist historians about what the Egyptians "must have seen" as a trivial, intuitive conclusion.

The $F_T$ formula has been praised for its "generality", "symmetry", "simplicity" and "exactness". Its "exactness", however, is based on the rather "coincidental" fact that the closely related $\frac{1}{3}$ – formula $F_P$ for the full pyramid is "exact," although it needs for its proof intrinsically "infinitesimal considerations" (see below, section 4). For the Egyptians to assume the factor $\frac{1}{3}$ may have been a matter of experience and convenience, not of deeper mathematical insight. But this makes the use of the word "exact" for the combined formula $F_T$ problematic as far as this use conjures up specific mathematical knowledge. As to the "generality" of $F_T$: this is – to some degree at least – a property of any formula which can really be called "mathematical" and thus it is in a way self-explanatory. "Symmetry" and "simplicity" of the formula for the truncated pyramid, however, can be disputed, as we will see below.

The current state of lack of knowledge about the way how the Egyptians derived the formula $F_T$ for the truncated pyramid is summarised in a recent book by a leading specialist:

> "There have been multiple attempts to explain how the Egyptian formula to calculate the volume of a truncated pyramid could have been achieved; some of these explanations use a clever modification of algebraic formulas—which were not used in Egyptian mathematics. Other explanations have tried to use practical experiences in form of woodworking. Another possibility would be the use of our knowledge of calculations of volumes of other objects from the hieratic mathematical texts. There are examples of calculations of the volumes of cylinders (e.g., papyrus Rhind, numbers 41 and 42) and a cube (papyrus Rhind, number 44). The underlying strategy for calculating their volumes is a multiplication of base and height. The procedure for calculating the volume of a truncated pyramid may be seen as a variation of that basic concept. Since the shape of the object indicates that a simple multiplication of the base and the height would result in a volume too big (i.e. that of the respective cuboid), a modification is put into place using three different 'bases' (the actual base (the square of the lower side), a rectangle made up of lower side and upper side, and the top area (a square of the upper side). To balance the sum of three areas (instead of one), they are then multiplied not by the height but by a third of the height." (Imhausen 2016, 75/76)

First, Imhausen alludes critically (as I will do below) to the work of historians around 1930, who ahistorically used a "clever modification of algebraic formulas" for deriving $F_T$ although such formulas were not used in Egyptian mathematics. Second, she allows the possibility that the formula originated from manipulating with models of parts of the frustum ("woodworking"), an interpretation which I will support below as well. Thirdly, Imhausen supports as "another possibility" the interpretation of the formula ($F_T$) as expressing an average of three volumes, or as being calculated via an average of areas multiplied by a common height.[10] This interpretation is very "intuitive" in the sense of considering the maximum and minimum volumes or areas which could be intuitively used for the calculation of the volume of the frustum and "balancing" it them by a third "intermediary" volume or area between the two extremes. However, this "intuitiveness" is not strong or complete enough suggesting the concrete value of the intermediates (*abh* or *ab*). Although, *ab* is the "geometrical mean" of the outer and inner squares $b^2$ and $a^2$, it remains mysterious why this value, taken as a "balance", delivers the exact volume. Still I will argue below that this kind of "inexact or approximate intuitiveness" was probably the main reason that the Egyptians preferred this formula over an "alternative" one which would result from directly dissecting the frustum.

One thing is the *interpretation* of a given formula (for instance as area- or as volume-average) and to find a hypothetical proof (within a given historical context) showing that this interpretation solves the original problem (here: calculation of a volume). From a pedagogical standpoint it is equally interesting although probably even more difficult to understand how the formula could have been found in the first place. I will reflect on this in the final chapter of this paper. Since an existing formula usually provides much inspiration and information to find a proof[11] and since a proof does not necessarily coincide with the original heuristics which led to the formula, it is pedagogically useful to hypothetically discuss the interplay of the four different cognitive constituents, namely problem, heuristics (discovery), theorem (formula), and proof, at a simple historical example like this.

[10] Things are a bit different in the Chinese case where the formula in the "Nine Chapters" points more clearly to averages of volumes not of areas (see below).

[11] As I argue below this information was insufficiently used by some historians of Egyptian mathematics around 1930 but better used by the Chinese commentator on the Nine Chapters, Liu Hui.

The lack of additional documents from the Egyptians has one advantage: it opens up for certain hypotheses which might have didactical value, assuming that modern pupils and ancient Egyptians were in a somewhat similar position with respect to learning mathematics.

Throughout my paper I am above all interested in tracing or conjecturing what I would call "germs of mathematical rationality", not however in finding full proofs in the modern sense. As I understand it this is the aim of some professional historians of ancient mathematical cultures as well (Chemla 2012). Anyway, we have but rarely the opportunity, or time, to present full proofs to our school pupils. And as teachers we should try to connect to the pupils' intuitions in the first place. We have many more historical documents from later historical periods such as Chinese mathematics in the first centuries of the CE, that give us more detailed information or at least reason to conjectures about mathematical intuitions in ancient people who were working on similar problems and found similar formulas (Wagner 1979, Chemla 2012). Nevertheless, Egyptian pyramids continue to fascinate the modern mind – maybe partly due to the lack of documents for this period – and the formula for the truncated pyramid retains some exceptional didactic flair.[12]

Before entering the main discussion in sections 5 through 8, I will remind the reader in section 3 of the concrete historical occurrence of the Egyptian formula $F_T$ and in section 4 of the need of "infinitesimal considerations" in order to derive the $\frac{1}{3}$ – formula $F_P$ mathematically. In the first case I am basing myself on known facts and literature, in the second case I reproduce and comment on Otto Neugebauer's opinions on this matter.

[12] The master pedagogue and mathematician George Pólya, for one, recognized the intuitive potential of the historical problems of the full and of the truncated pyramid and commented on them in *Mathematics and Plausible Reasoning* (Pólya 1968, vol.2, 162, 188) as well as in *Mathematical Discovery* (Pólya 1981, vol. 2, 2-10). But he used modern, albeit simple, algebraic methods, and did not attempt to make a connection with a possible Egyptian way of deriving the two formulas. In (Pólya 1981) this is part of a rather complicated theory of heuristics whose usefulness for the classroom I doubt. A historical view would have helped.

## 3. The "formula" in the Moscow papyrus

The specialists among historians of mathematics largely agree that problem no. 14 in the so-called Moscow papyrus shows the familiarity of the Egyptians with an algorithm (formula) for calculating the volume of a truncated Egyptian (quadratic) pyramid.

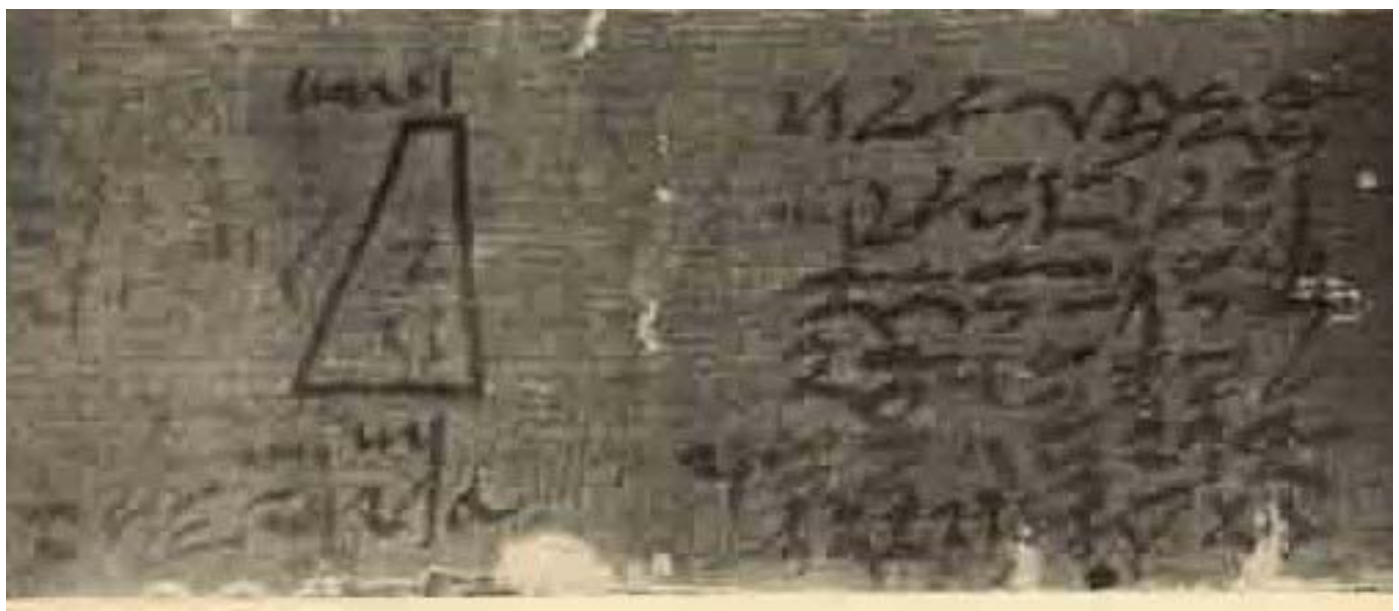

Calculation of Volume of Truncated Pyramid.
XIIth dynasty.

**Figure 3:** Picture from the first (English) publication on problem 14 of the Moscow papyrus (Touraeff 1917, 101)

In his book of 1934 "Pre-Greek mathematics" (*Vorgriechische Mathematik*), which has apparently never been translated into English, the Austrian Otto Neugebauer, the most famous historian of ancient mathematics and astronomy in the 20th century,[13] says:

> "The jewel in the crown of Egyptian mathematics is the correct formula for the volume of a truncated pyramid with a quadratic top area." (Neugebauer 1934, 126)

In algebraic notation (replacing the actual numbers used in the papyrus) the verbal description in the papyrus – which is reproduced further below – corresponds to the above formula $F_T$.

In his 1934 book Neugebauer also explains what he meant by "formula" for the historical period under discussion, i.e., in the beginning of the second millennium BCE:

> "When we speak of a 'formula', then this is of course to be understood here as well as otherwise always in such a way that the text itself calculates in actual numbers, but according to a rule which is just expressed by our formula." (Neugebauer 1934, 127)

[13] Neugebauer began as a historian of Egyptian mathematics in the 1920s, but he became famous as the interpreter of the Babylonian astronomical and mathematical cuneiform tablets. For the most recent evaluation of his life and work see (Jones/Proust/Steele 2016).

Although the figure in the Moscow papyrus at problem 14 is not drawn symmetrically,[14] historians assume that it refers either to the symmetrical truncated pyramid (Figure 1 above) or to the special asymmetrical truncated pyramid with two perpendicular triangular faces. The latter could be called a "right" asymmetrical pyramid; four of them form a symmetrical frustum with the same height $h$ and a double base $2a$. This provides an additional intuitive argument for the fact that the "right asymmetrical pyramid" has the same volume as the symmetric one, because doubling the base suggests quadrupling the volume for the same height. For this conclusion Cavalieri's principle is needed in some simple form considering the volume as composed of square areas. Below I will argue – based on a suggestion by Shutler – that even without use of Cavalieri's principle a very simple intuition of similarity for volumes works to establish the $\frac{1}{3}$ – formula $F_P$ for symmetric pyramids and following from this also for $F_T$.

The main reason for the claim that the Moscow papyrus contains such a formula $F_T$ with general meaning is of course the elaborateness of the numerical algorithm (see below), which can easily be applied to pyramids of any dimensions provided $a > b$. The numbers for the units used in the text $a = 4$, $b = 2$ and $h = 6$, may look strange, because they belong rather to a truncated obelisk than to an Egyptian pyramid as we know it from Giza.[15] But perhaps these extreme numbers only confirm the impression that the formula is supposed to have general meaning.[16] Or the numbers were chosen because they were the simplest (natural) numbers, which produce a simple result of calculation.[17]

According to the Egyptian text, the algorithm proceeds as follows, starting with $a^2 = 16$, $ab = 8$ and $b^2 = 4$, where the values within the parenthesis in the above formula $F_T$ are summed first. Regardless of whether one considers the formula as the mean value of volumes

---

[14] It has been pointed out by detailed research that diagrams in old mathematical manuscripts are not always in accordance with the accompanying mathematical text. Ken Saito, Nathan Sidoli and others have shown this for early documents connected to Euclid's *Elements* (Saito/Sidoli 2012). Karine Chemla reminds me that these figures are often more symmetrical than the text suggests. Symmetry in picture implies restriction of generality while in the case of the asymmetrical frustum the figure may indicate generality and possibly understanding of Cavalieri's principle. De Young (2009, 350/51) considers also the diagram for problem 14 of the Moscow Papyrus but does not comment on the lack of symmetry there.

[15] (De Young 2009, 351) refers to the extreme steepness of the diagram.

[16] Shutler (2009) draws a special conclusion from the fact that $b$ is half of $a$ in the example. See below. The "generality" of the formula includes that – similar to nearly the same formula in the Chinese Nine Chapters – the corner pyramids of the frustum are not inscribed in cubes but more generally in cuboids (see below).

[17] This is a conjecture in (Gunn/Peet 1929, 184/85). The choice of 6 instead of 3 as the height, which leads to the "atypical" shape of the frustum, might be related to the wish to distinguish this number from the divisor 3 which belongs as a constant to the formula. This was suggested to me by my student Gunnar Andersen. However, we have a similar situation in Papyrus Rhind for the calculations of circular areas where the number 9 is used for two different purposes, the diameter and the factor in a formula. Since there is another example with diameter 10 in this papyrus it is nevertheless relatively easy there to recognize the "formula".

(an opinion I tend to share) or as the mean value of areas which is afterwards multiplied by a common height $h$, the formula is based on a certain elementary "algebraic thinking" of the Egyptians, by taking the factor $h = 6$ from the sum of volumes, and by applying the average $\frac{1}{3}$ not directly to the sum of areas but first to the factor $h$:

(1) Add 16 plus
(2) 8 and in addition 4
(3) so you get 28.
(4) Calculate $\frac{1}{3}$ of 6. So you get 2
(5) Calculate two times 28. So you get 56.
(6) See: the volume is 56. You have calculated correctly!

According to the literature[18] the numbers relate to cubits (ca. 457 mm), thus the height $h$ of the frustum is about 3 m for the example from the Moscow Papyrus. This can be hardly considered to refer to a model. The experimental verification of the formula would require smaller dimensions, such as in the Chinese case (see below section 6).

Historians agree that there is no evidence of a plausibility argument, let alone a proof of the $F_T$ formula in the Moscow Papyrus or elsewhere in Egyptian mathematics. So if we speak of an "algorithm" here we do not mean an algorithm with explanatory power in a sense of a geometrical construction, as we find it in Euclid and in Liu Hui much later.

## 4. "Infinitesimal considerations" needed to prove $F_P$

In a 1933 article[19] in his own journal *Sources and Studies for the History of Mathematics* (Quellen und Studien zur Geschichte der Mathematik), in which Wasili Struve had published the first German transcription of Problem 14 of the Moscow Papyrus three years earlier

[18] In (Struve 1930, 135) cubits (Ellen) are assumed. But this is followed by the explicit remark (p.137) that the Moscow Papyrus does not give names for the length units.
[19] The paper (Neugebauer 1933) does not focus on Egyptian mathematics but rather on a Babylonian formula for the volume of the truncated pyramid. So his remarks about the Moscow Papyrus should be understood as side-remarks and thus not be unduly criticized. Nevertheless, Neugebauer was, particularly in the first period of his career, deeply and creatively involved also in the philological interpretation of Egyptian mathematics texts, as repeatedly and gratefully acknowledged in Struve (1930).

(Struve 1930), Neugebauer is somewhat condescending about the many attempts by historians of his time to find out how the Egyptians might have derived the formula $F_T$:

> "Much effort has been made to give a derivation for the exact formula, an effort that does not seem particularly interesting to me, since it basically amounts to developing a fundamentally wrong 'proof', without being able to show why precisely this and no other error was made in the derivation." (Neugebauer 1933, 347/48)

I do not particularly like this statement by Neugebauer for several reasons. It seems to be based on the "presentist" assumption of a timeless "exactness" in mathematics and seems to proclaim a ban on hypotheses, although speculations and suppositions can be very helpful for historical as well as for didactical purposes. Moreover, Neugebauer is not consistent, for he himself starts speculating immediately after this statement, as we shall see below. Most importantly, Neugebauer does not clearly distinguish between the two problems connected to the two formulas and thus discourages from the outset efforts to prove $F_T$ under the assumption that $F_P$ is correct.

In his 1934 book Neugebauer continues to discuss the two different problems connected to $F_T$ and $F_P$ although ostensibly he is talking only about the former:

> "What is surprising about this formula is above all two things: on the one hand the symmetrical shape, on the other hand the mathematical correctness, which, as is well known, requires infinitesimal considerations, i.e. goes beyond the framework of elementary geometry, especially with this formula, if it should be derived correctly." (Neugebauer 1934, 126)

The second point, the "infinitesimal considerations", is obviously related to the more special formula $F_P$. "Infinitesimal" can here only be understood in a very broad sense as involving "infinitary" arguments or limits, maybe by indirect proof, not however, claiming a rigorous theory of "infinitesimal analysis" in the modern sense. Neugebauer must have included Euclid here because he was of course aware, that Euclid in the 3rd century BCE in book XII of his *Elements* had given a perfectly rigorous proof (according even to the modern standard of Neugebauer's time) of the $\frac{1}{3}$ – formula $F_P$. But Neugebauer's implicit assumption that Euclid's proof would not have been in reach of Egyptian mathematics seems justified.[20]

Moreover, it seems certain that Neugebauer's claim about the requirement of "infinitesimal considerations" reflects, in addition, even more modern mathematical results of his own time of which he was – as a mathematician trained in Göttingen – well aware. Above all one has to think of Max Dehn's solution (1900) of Hilbert's third Paris problem of the same year, which shows that two arbitrary polyhedra (in particular tetrahedra two of which make a pyramid) of equal volume cannot as a rule be dissected in finitely many steps into pairwise congruent parts.[21] Therefore one cannot expect that arbitrary pyramids can be shown by dissection to be equal in volume with simple figures such as cuboids. This does not rule out that *particular* pyramids (for instance Juel pyramids) can be shown to be parts of simple figures which can be calculated.

Neugebauer, although ostensibly arguing against speculations, says in his 1934 book rather vaguely:[22]

> "The most likely assumption seems to be that they found for a particular and simple case (for instance the corner pyramid in a cube, see fig. 41) by an intuitive argument the right formula and then extended it to the non-trivial case." (Neugebauer 1934, 128)

In the accompanying figure 41 of his book which we call here figure 4, Neugebauer shows the dissection of a cube into three congruent asymmetrical pyramids:

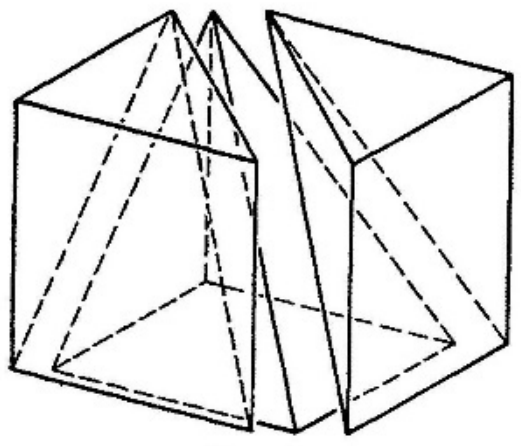

Fig. 41.

**Figure 4:** A cube dissected into three "right" asymmetrical pyramids (Neugebauer 1934, 128)

[20] Euclid used in the typical Greek manner indirect proofs and exhaustion methods coming from Eudoxus. Euclid's proof also uses the pairwise equality of the volumes of three tetrahedra with pairwise equal bases and heights, from which he concludes the general equality of the three. This would probably have put even greater demands on the Egyptians' abstract mathematical sense than the more intuitive Cavalieri principle. Gunn and Peet conclude that the Egyptians would not have had an understanding for such typically "modern" mathematical abstraction as in Euclid XII, 5/7 (Gunn/Peet 1929, 181).

[21] This is different from the analogous problem for areas in the plane, as results from the 19th century show.

[22] Probably Neugebauer means by "extension to the none-trivial case" something like "stretching" a cuboid and the pyramid in it by a linear factor.

Unfortunately, the pyramids built in Egypt did not have the same height $h$ as the side $a$ of the square base $a^2$. Therefore the "trivial case" had to be extended. The Egyptian pyramids at Giza are approximately halves of regular octahedra.[23] For general straight cuboids, which have two squares and four rectangles as faces, a comparatively simple dissection as in the case of the cube cannot be expected.[24]

As to the problem to derive $F_T$ under the assumption of the validity of $F_P$ the situation is of course different and Euclid's indirect methods or any other "infinitesimal considerations" are not needed in this case.

In this sense, claiming an independent historical importance for $F_T$, Touraeff's remark of 1917, in the first (English) publication of the Moscow Papyrus, also seems justified:

> "If only our explanation of the problem is right, we have here a new and interesting fact, i.e. the presence in Egyptian mathematics of a problem that is not to be found in Euclid." (Touraeff 1917, 102)[25]

## 5. Neugebauer, Gunn and Peet, and the dangers of a "presentist" misrepresentation of the truncated pyramid formula

In his 1934 book Neugebauer presents a dissection of the asymmetrical truncated pyramid by reasoning which he apparently considers to be intuitively plausible and in reach of Egyptian mathematicians:

> "The conditions on which this consideration is based are such that they can be expected of Egyptian mathematics. The geometrical dissection is immediately evident in the case of asymmetrical bodies." (Neugebauer 1934, 128)

---

[23] The measurements for the Cheops pyramid are about $a$ = 230 m and $h$ = 145 m.

[24] It was the BBC's video "The story of mathematics" which I occasionally use in my lectures that stimulated my interest in the pyramid problem. Here it is implicitly and erroneously claimed that any cuboid can be cut into three congruent asymmetrical pyramids, which allegedly could have been a possible solution for the Egyptians.

[25] We might add that Euclid like the Egyptians had no algebraic methods as we have them now. So it would not have been trivial though not difficult for him to derive $F_T$ from $F_P$. As is well known, a systematic translation of geometrical statements into algebraic, symbolic language was only introduced by Descartes in the 17th century.

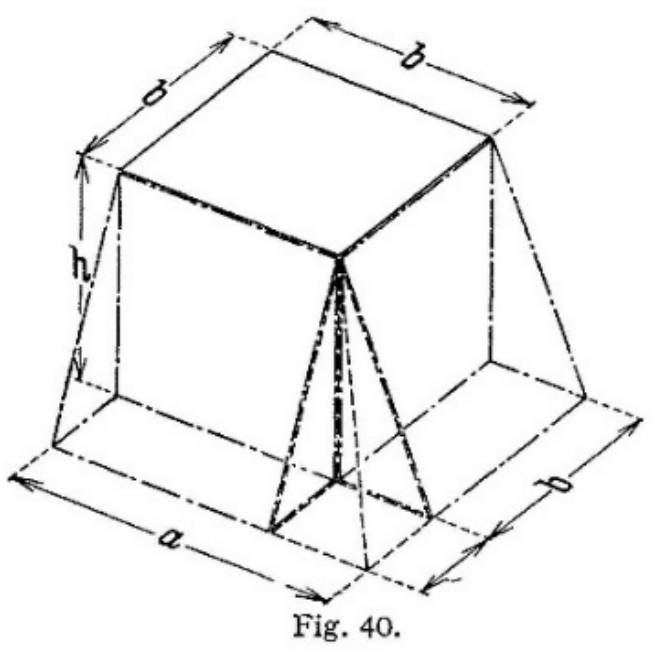

**Figure 5:** from Neugebauer (1934, 128) for the "right" asymmetrical truncated pyramid

What Neugebauer means is the choice of the frustum of the "right" (as I called it above) asymmetrical pyramid in figure 5 as the starting point and putting one of the two triangular prisms with volume $b\,(a\text{-}b)\,\frac{h}{2}$ on top of the other which results in a box[26] with the volume $hab$. There is in addition one small asymmetrical corner pyramid with base $(a\text{-}b)^2$. Assuming the knowledge of $F_P$ for the latter one has (in algebraic notation which for our convenience stands here for a purely geometrical dissection) the following expression for the volume of the asymmetrical frustum:

$$V_T = h\,b^2 + 2\,b\,(a\text{-}b)\,\frac{h}{2} + \frac{h}{3}\,(a\text{-}b)^2$$

$$= hab + \frac{h}{3}\,(a\text{-}b)^2 \qquad (F_{TA})$$

(henceforth $F_{TA}$ for "truncated-alternative")

In contrast to Neugebauer I believe that the geometrical dissection of the symmetrical truncated pyramid of figure 1 is just as obvious and simple as the one for the right asymmetrical frustum. The dissection of the symmetrical frustum was indeed assumed in 1929 by Gunn and Peet as a way for the Egyptians (see below). The four asymmetrical corner pyramids of the symmetrical truncated pyramid (figure 1), each of which having the volume $\frac{h}{3}\left(\frac{a-b}{2}\right)^2$, can be easily assembled to form a pyramid on the square with the area $(a\text{-}b)^2$. This requires no algebraic manipulation because these volumes are represented by simple geometrical figures. And the resulting small pyramid is even symmetrical while the resulting formula $F_{TA}$ for the volume is identical to the one given by Neugebauer. The latter author seems to have assumed Cavalieri's principle to be intuitively convincing for the Egyptians because he finds the case of "asymmetrical bodies" as unproblematic and "immediately evident". This in itself seems to me highly "speculative".

[26] "Box" is in the following used for a further generalization of "cuboid", and is basically synonymous with parallelepiped, which does not necessarily have any quadratic faces.

None of the historians of mathematics around 1930 was of course so naïve as to assume that one can now go on and simply replace the squared binomial $(a\text{-}b)^2$ in $F_{TA}$ by the formula $a^2 - 2ab + b^2$. From this the volume formula $F_T$ would just follow. Historians knew that there was no algebra in our modern sense in Egypt at the time and that the formula $(a\text{-}b)^2$ would have to be interpreted geometrically in the plane. But such an interpretation would be highly speculative in the case where the plane areas did not have any intuitive connection to the frustum. In addition, it was not quite so easy to find a geometric alternative form of $(a\text{-}b)^2$ which is in the formula $F_{TA}$ as it would have been the case for, let us say $(a+b)^2$.

In his 1933 article, immediately after his refusal to make hypotheses, Neugebauer tries nevertheless to speculate about the algebraic character of the original formula $F_T$:

> "Important seems… the symmetrical shape of the formula [$F_T$], which looks like an algebraic transformation of an initially differently constructed expression, without it being possible to see how something like this should be achievable with the computational aids known to us from other Egyptian mathematical texts." (Neugebauer 1933, 348)

A year later, however, in his 1934 book, Neugebauer comes with additional information on "other Egyptian mathematical texts" which he seems to have found or interpreted since then:

> "The only transformation that is necessary is the squaring of the binomial $a\text{-}b$, and it is exactly this that is documented in other [Egyptian] texts." (Neugebauer 1934, 128)

Using these "other texts", which were apparently unrelated to the calculation of the volume of pyramids,[27] Neugebauer was able to express the volume of the corner pyramid in the following way

$$\frac{h}{3}\,(a\text{-}b)^2 = h\left(\frac{a^2}{3} - \frac{2}{3}ab + \frac{b^2}{3}\right)$$

Neugebauer does not try to replace his algebraic argument by more intuitive geometrical reasoning. According to $F_{TA}$ adding the volume of the box *hab* results in the formula $F_T$ of the Moscow Papyrus.

[27] Neugebauer refers here to an earlier passage in his book (p.112), which however does not reveal which Egyptian text he is referring to.

Neugebauer insinuates that the symmetry in $F_T$ comes from an "algebraic transformation" of a part of formula $F_{TA}$. However, I find that the formula $F_T$ in its entirety can be more intuitively and pedagogically more usefully interpreted as an average of either volumes or areas (multiplied by heights). If one (like Neugebauer) assumes the possibility or even the necessity of "algebraic" arguments, I think that somehow the historical flair and educational potential of a discussion about $F_T$ is destroyed. Anyway, it seems to me that one should historically always look first for the simplest, most intuitive explanation, which presupposes the least formalized mathematics.

In a vein similar to Neugebauer, the British Egyptologists Battiscombe Gunn and Thomas E. Peet proceeded in an earlier article of 1929, but without having Neugebauer's additional information about "other texts". Their article (Gunn/Peet 1929) seems to me the most substantial and illuminating discussion of the historical formula for the truncated pyramid before (Wagner 1979), the latter, however, relating to a totally different historical context.

The two Brits consider different dissections of the symmetrical truncated pyramid. They assume that the Egyptians had to resort to using pyramid-shaped models made of dried mud to find the volume of the pyramidal parts at the corners of the figure:

> "By a method remote from those of pure geometry, let us say by cutting up a lump of half-dry Nile mud with a piece of stout thread (180), … we may suppose that he [the Egyptian] had recourse to weighing." (Gunn/Peet 1929, 180/182)

Interestingly but not really surprisingly, Gunn and Peet end up mathematically with the same formula $F_{TA}$ for the truncated pyramid as Neugebauer would five years later in his book (see above) and which Liu Hui had obtained around 250 CE as the second of two alternatives (see below section 6):

$$V_T = h\,a\,b + \frac{h}{3}\,(a-b)^2$$

In contrast to Neugebauer, however, the British authors are more articulate or imaginative about how the second term of the equation can be found, namely by weighing and equating the volume of the pyramid to one third of a small cuboid on the same basis.

But – again like Neugebauer and although they stressed the use of models heavily – the British authors went one step too far for being able to derive the formula of the Moscow

Papyrus. By combining the elementary parts of the frustum into figures (parallelepiped and combined corner pyramid) they restricted the possibility of further manipulation with them. Instead, they assumed that the Egyptians would now have seen the need to “simplify” $F_{TA}$ by taking out the common factor $h$ which leads to

$$V_T = h\,[a\,b + \frac{1}{3}\,(a-b)^2]$$

This elementary arithmetic simplification for avoiding repeated multiplication which is reminiscent of the algorithm in the Moscow Papyrus (section 3 above) was according to the authors “well known to the Egyptians” (Gunn/Peet 1929, 182), also documented in other texts as in the Rhind Papyrus. This cannot be disputed. However, applying this operation at this point restricted even more the possibility of reaching formula $F_T$, because now every next step of “simplification” would have to proceed on the level of manipulations with the area inside the parenthesis, instead of with the much more intuitive elementary three-dimensional figures derived from the frustum.[28] Above all, however, the classical formula $F_T$ is much more “intuitive” in the sense of relating to averages of areas or volumes, as indicated in the quote from Imhausen above. By way of contrast the area $a\,b + \frac{1}{3}\,(a-b)^2$ in the “simplified” form of $F_{TA}$ has no relation to the intuitive maximum or minimum of areas for the base of the frustum which the average $\frac{1}{3}\,(a^2 + ab + b^2)$ with the same numerical value has, which occurs in $F_T$.

But the formula $F_{TA}$ is very simple, and the Egyptians could have stopped here. Indeed, the British authors say:

> “At this point the Egyptian [sic], whose solutions of problems were often very cumbrous might well have rested content.” (Gunn/Peet 1929, 182)

However, the authors then give various reasons why this should not have been the case, why the Egyptians should have further “simplified” the expression to finally produce the formula $F_T$. As the main reason the authors name educational concerns:

> “The present writers suggest that this elementary simplicity [of the formula $F_T$; RS] was given to the problems that they might be the easier to learn by heart; that the education of the student consisted partly in his committing them to memory in order to be able to apply them *mutatis mutandis*, to similar ones that might arise.” (Gunn/Peet 1929, 185)

[28] The more intuitive but arithmetically less efficient form of $F_{TA}$ as a sum of volumes is chosen by Liu Hui (Chemla 1991, 87).

My main question now would be: is this last expression $F_{TA}$ for the volume of the truncated pyramid – with or without taking out $h$ as a common factor – really more complicated than the formula $F_T$ in the Moscow Papyrus? It has only two terms, unlike the three in the papyrus, both formulas $F_T$ and $F_{TA}$ require that you take the square of a length and divide some quantity by 3. To find the difference ($a$-$b$), which is a positive geometrical length, should have been within the reach of the Egyptians and does not mean manipulating with negative numbers.

Both formulas $F_T$ and $F_{TA}$ are equally general and exact and allow an easy calculation of the volume of the full pyramid by setting $b = 0$. Interpreted as step by step algorithms none of the two formulas is "simpler" than the other. What is then the "elementary simplicity" which Gunn and Peet claim for $F_T$ ? First, the figures expressed in $F_T$ are all of a similar calculable kind, namely parallelepipeds of the same height $h$, making the formula in this sense more "homogeneous" when considered as a formula about three volumes. More importantly, as mentioned before, the classical formula $F_T$ is more "intuitive" as an average related to volumes or areas. This is a kind of "approximate intuitiveness," however, which does not allow conclusion that the value expressed in the formula is really correct.

It is unclear to me whether when they spoke about "easier to learn by heart", Gunn and Peet also hand in mind the "algebraic symmetry" of the papyrus formula, which is missing in the alternative formula $F_{TA}$ for the volume. The "ease of memorizing" that we associate with the formula $F_T$ may result precisely from its similarity with the expression for $(a+b)^2 = (a^2 + 2ab + b^2)$ with which we are so familiar. One historian (Thomas 1931) has tried – without much success – to take inspiration from the identity $(a–b)\ (a^2 + ab + b^2) = a^3 – b^3$.[29] This was another instance of the dominance of algebraic thinking (this time on the level of volumes) because it seems highly unlikely that the Egyptians would have found this formula by physical or mental manipulation with volumes. This would have required the fictious consideration (creation) and subsequent neglection (subtraction) of boxes such as $a^2b$ and $ab^2$ which are no immediate results from possible dissections of the truncated pyramid. The Chinese way (see below) is much closer to immediate spatial intuition.

---

[29] Interestingly enough, the volume of the *truncated* Juel pyramid can be represented as one sixth of the difference of the $a^3$ and $b^3$ cubes, the height $h$ being in this case $(a$-$b)/2$. So in this special case the identity $(a$-$b)$ $(a^2 + ab + b^2) = a^3 – b^3$ makes immediate geometrical sense for $F_T$. This is not mentioned in (Thomas 1931) or (Gillings 1964). Moreover, if one stretches both cubes and with it one of the six truncated Juel pyramids inscribed in $a^3$ with the same factor $t$, the height of the Juel pyramid increases by the same factor as well, and the general formula $F_T$ is easily confirmed. But all this still requires that the above identity is confirmed.

Perhaps the admiration for the "elementary simplicity" of the formula in the Moscow Papyrus is due to the fact that it is not (!) immediately recognizable as being connected to the dissection of a truncated pyramid, but rather reminds us of "modern" algebra? That we attach so much value to "symmetry" may just be a prejudice coming from our algebraic education.

In any case, Gunn and Peet (1929), suggesting that the Egyptians would *not* have been satisfied by the above formula $F_{TA}$ for the frustum, looked out for an Egyptian way to simplify the hypothetical formula. As already mentioned, they first took the common factor $h$ out of the sum, an operation which in their opinion was within the reach of the Egyptians. Then they resorted to careful and lengthy considerations of an algebraic and experimental kind but related to areas, not to volumes. They came to the conclusion that three rectangles and a square, drawn together on papyrus paper – by cutting and recombination – might have suggested the following identity, which immediately leads to $F_T$:

$$3\, ab + (a-b)^2 = (a^2 + ab + b^2)$$

This speculation, which is reminiscent of an algebraic interpretation of some geometric identities in Euclid's *Elements*, shows once again that some historians of Egyptian mathematics around 1930 were not able to distance themselves from geometric-algebraic manipulations familiar to them from the historiographic presentation of Euclid[30] when looking for the way to derive formula $F_T$. I believe that already the decision of Gunn and Peet – and later also Neugebauer – to write the experimentally (by weighing) determined value for the volume of the four corner pyramids in algebraic language as the second term $\frac{h}{3}\,(a-b)^2$ and not, for example, in the form of an icon reminding of a pyramid has led and seduced them into algebraic thinking.

*To sum it up in a more principled manner:*

In my opinion, most historians who originally, around 1930, commented on the Moscow Papyrus (Struve, Neugebauer, Vogel, Thomas, Gunn/Peet)[31] and who discussed possibilities for an original Egyptian proof of formula $F_T$, were influenced in their historical judgment by some "presentist" thinking which they shared. To me there is a remarkable contrast between the philological depth and shrewdness of these historians and the relative lack of methodological sophistication and reflection when it came to mathematical interpretation. I found this most strikingly exemplified by Struve (1930), although I will not go into his mathematical commentary in detail which seems to me much less sophisticated than those by Gunn, Peet and Neugebauer. Of course this is for the most part just an expression of the relative youth of the historiography of Egyptian mathematics compared to Egyptology as a whole.

[30] This is the disputed interpretation of Euclid's geometry as a kind of "geometric algebra".
[31] I cannot discuss all of them here in detail.

The word "presentist" has hitherto been used in this article in a rather loose manner. It shall therefore now be explained in somewhat more detail. First of all it does not primarily mean judging the past from the (then) "present" of the 20th century or even insinuating the existence in Egypt of algebraic symbolism; the latter was basically a product of the 17th century. It does not even mean the presence of knowledge of Greek mathematics from the 3rd century BCE, at least not of the more sophisticated Greek methods such as Eudoxus's proportions and exhaustion. However, to claim that those "historians around 1930" insinuated Greek geometry comes closer to the meaning of "presentist" thinking which I have in mind. Those historians tended to write down assumed Egyptian formulas in symbolic algebraic notation (which is not to meant as criticism because it was the language they were used to) but then interpreted these formulas in terms of their idea of Euclid's plane geometry without really caring anymore about connection to the original three-dimensional problem of the frustum. In this sense the "presentist" approach had also something "algebraicist" to it, prioritizing manipulation of formulas over geometric understanding. Part of the problem was probably the small "algebraic" step which is actually contained in the formula of the Moscow Papyrus which puts the height *h* at the end of the algorithm (and has thus to be written first in the formula) for the convenience of the calculation (as noted above also in connection with the alternative formula $F_{TA}$).

By way of contrast to such a "presentist" interpretation, I will look in the following section 6 for another, alternative, simpler and more intuitive derivation of the volume formula $F_T$ for the truncated pyramid, which in my opinion would strictly fall within the scope of the possibilities of Egyptian mathematics, provided that $F_P$ is known.

In section 8 I will then – after some methodological reflection in section 7 – reiterate Shutler's recent (2009) intuitive proof of $F_P$, which in my opinion would not have been out of reach for the Egyptians either.

In both cases, with respect to $F_T$ and $F_P$, I do not reject some mild kind of speculations, encouraged by the fact that specialists like Neugebauer, Gunn and Peet have speculated constantly as well.

## 6. Returning to the original volume problem for the truncated pyramid in the manner of Liu Hui (263 CE)[32]

In this section I will discuss – in broad outline and quoting mostly indirectly from recent scholarly work on Chinese mathematics – Liu Hui's proof and try to show its potential for explaining the formula in the Moscow Papyrus. While thus taking inspiration from Liu Hui, I will at the same time connect to several points of my interpretation in the preceding section 5, such as the juxtaposition of the classical (documented) formula $F_T$ and the "alternative" formula $F_{TA}$. I find exactly in Liu Hui's commentary such an alternative formula which is nearly identical with the one suggested by the historians of Egyptian mathematics around 1930.

The classical work "Nine Chapters" of Chinese mathematics was compiled in its most authoritative form as a manuscript written by unnamed Chinese mathematicians in the first century CE. Throughout history the Nine Chapters have been commented upon both by the Chinese themselves and most recently in Western or bilingual editions. For all the nuances which cannot be discussed here I refer to the monumental French bilingual edition (Chemla/Guo 2004) which comprises 1141 pages, and to many articles by Karine Chemla. The much shorter English-Chinese edition (Guo/Dauben/Xu 2013) is particularly useful for English readers. For quotes I will rely here on (Wagner 1979) which is best accessible and is based on his pioneering Copenhagen master thesis (Wagner 1975). The modern editions, both in Western languages or in Chinese usually contain not only the classical text of the Nine Chapters but also historical Chinese commentaries on this text of which the commentary by Liu Hui (around 263 CE) is most famous. While the classical text of the Nine Chapters as a rule formulates problems and provides algorithms for their solution without proof, Liu Hui always gives explanations and adds own problems and solutions.

In the Nine Chapters the formula for the volume of the truncated pyramid appears in chapter 5 as solution to problem 10 and can be written as

$$V_T = \frac{1}{3} \bullet h \bullet (a^2 + ab + b^2)$$

with the dot • denoting multiplication.

[32] After reading an earlier version of this paper my colleague Reinhard Bölling (Berlin) convinced me to draw the parallel to Liu Hui more clearly.

Like in the Egyptian case there is only given one numerical example. However, the algorithm in its detailedness is again convincing enough to write it down as the above formula. The only difference compared to the Egyptian $F_T$ formula in the Moscow Papyrus is that the factor $\frac{1}{3}$ is applied separately at the very end of the algorithm (therefore written at the beginning in our notation according to the usual understanding of multiplication). This is no minor point because it makes its interpretation as an average of volumes very likely from the outset.

Indeed, the main point of Liu Hui's commentary is that, unlike the historians of Egyptian mathematics around 1930, he took the historical algorithm in the Nine Chapters (i.e. the formula) really seriously in the sense that he associated it with the original volume problem of the frustum. As I understand it, he read it as a sum of three volumes divided by 3 (i.e. an average of three volumes) and wanted to relate these volumes to the volume of the truncated pyramid and thus considered three copies of the frustum. What the original algorithm in the Nine Chapters itself does not do is explaining how the parts of the frustum with their concrete measures could be related "one to one" to the volumes of the simple boxes in the formula. By revealing this relationship Liu Hui confirmed the concrete numerical value of the volume as explicitly given in the Nine Chapters as the result of the formula. In other words: "algorithms" for calculating the volume of the truncated pyramids are given both in the Moscow papyrus and in the classical "Nine Chapters," and the formulas are nearly identical. But a "derivation", i.e. an explanation is given only in Liu Hui's commentary.

Starting point for Liu Hui was the formula in the Nine Chapters (above) and the concrete example given there was $a = 5$, $b = 4$ and $h = 5$. Here the unit length was "zhang" which is 10 "chi". The latter can be assumed to have been about 23 cm for the third century CE, according to (Chemla/Guo 2004). This amounts to a height of 11.50 m and a total volume of 101,666 + 2/3 cubic chi (this result is explicitly given in the Nine Chapters) which corresponds today 1236.97 $m^3$. These huge dimensions in the original "Nine Chapters" were probably chosen in order to demonstrate that one really needs a formula for calculation of the volume because a recourse to dissecting or weighing was impossible. The other way round, to make the formula plausible with the help of models which can be handled[33] one had to reduce the dimensions conspicuously, and this is what Liu Hui did. The major problem were the extremely steep corner pyramids in the original example with $(a\text{-}b)$ 2 = 0.5 zhang as base edge and $h = 5$ zhang as height. These dimensions were far from a pyramid inscribed in a cube or from a Juel pyramid with $45^0$ slope of the faces, and these were the two cases of pyramids where the volumes can be easily calculated. Therefore Liu Hui chose not only smaller dimensions but such dimensions which produced simple corner pyramids.

[33] "Handle" does here not necessarily mean physical manipulation with models, although Wagner assumes such handling: "It is clear that he owned a set of these blocks himself and expected his readers to use them also." (Wagner 1975, 13) This interpretation has been doubted in more recent commentary.

Liu Hui therefore based his explanation on a very simple example[34] where all nine "parts" of the frustum were so-called "standard blocks" with different shapes (not to be confused with "boxes" in the sense we discussed before) with height 1 and base 1, namely cubes, prisms (qiandu) and right asymmetrical pyramids (yangma). This corresponds to the example $a = 3$, $b = 1$ and $h = 1$, in each case assuming "chi" ≈ 23 cm as the length unit.

Liu Hui produced the standard blocks, which were also used for explaining other volume problems than the truncated pyramid, by first cutting a cube through the diagonals of opposite quadratic faces into two (special) prisms and then cutting from one of the prisms two thirds, which results in the right asymmetrical pyramid (yangma).[35] In the given simple example the corner pyramids are inscribed into cubes and can therefore (at the end of the procedure to be described now) be easily integrated into the three boxes of the formula $F_T$ of the Nine Chapters. The relation $(a-b)/2 = h$ [= 1] is essential for this property. By way of contrast the relation $b = h$ [=1] which makes the four prisms of the frustum to half cubes, and the ratio $a/b = 3$, which allows taking the base edge of the corner pyramids as unit are accidental to the volume problem. They could be replaced by any other ratios without impairing the general validity of Liu Hui's argument.

In order to establish the correspondence between the volumes of the three boxes in the formula $F_T$ of the Nine Chapters and the truncated pyramid Liu Hui starts with three copies of the latter with reduced dimensions in the sense of the standard blocks. He dissects one of the frustums into the nine parts mentioned above: one central cuboid, four prisms and four asymmetrical corner pyramids (Wagner 1979, 169). In this first step Liu Hui combines the central cuboid (which is here a cube) and the four prisms into the box *hab* thus producing the first of the three volumes in the formula $F_T$ (actually the one in the middle of it). This could be easily imagined mentally and corresponds fully to the dissection which led Neugebauer, Gunn and Peet to the "alternative" formula $F_{TA}$. However, unlike the latter formula, Liu Hui did not at this first step combine the four corner pyramids into a symmetrical pyramid with base edge (*a*-*b*).

[34] Karine Chemla prefers the use of "paradigm" for it.
[35] (Wagner 1979, 175).

Now there could have been many ways – at least for the particular and simple dimensions of Liu Hui's standard blocks – to combine the parts of the two additional copies of the truncated pyramid, together with the four right pyramids (yangma) which are left from the first frustum, into the other two volumes (boxes) of the classical formula $F_T$ which was of course always on Liu Hui's mind.[36] It was also clear from the classical formula that for Liu Hui's special example the volume was equal to 13 cubes divided by three, and Liu Hui mentions the corresponding volumes of resulting boxes at each step.

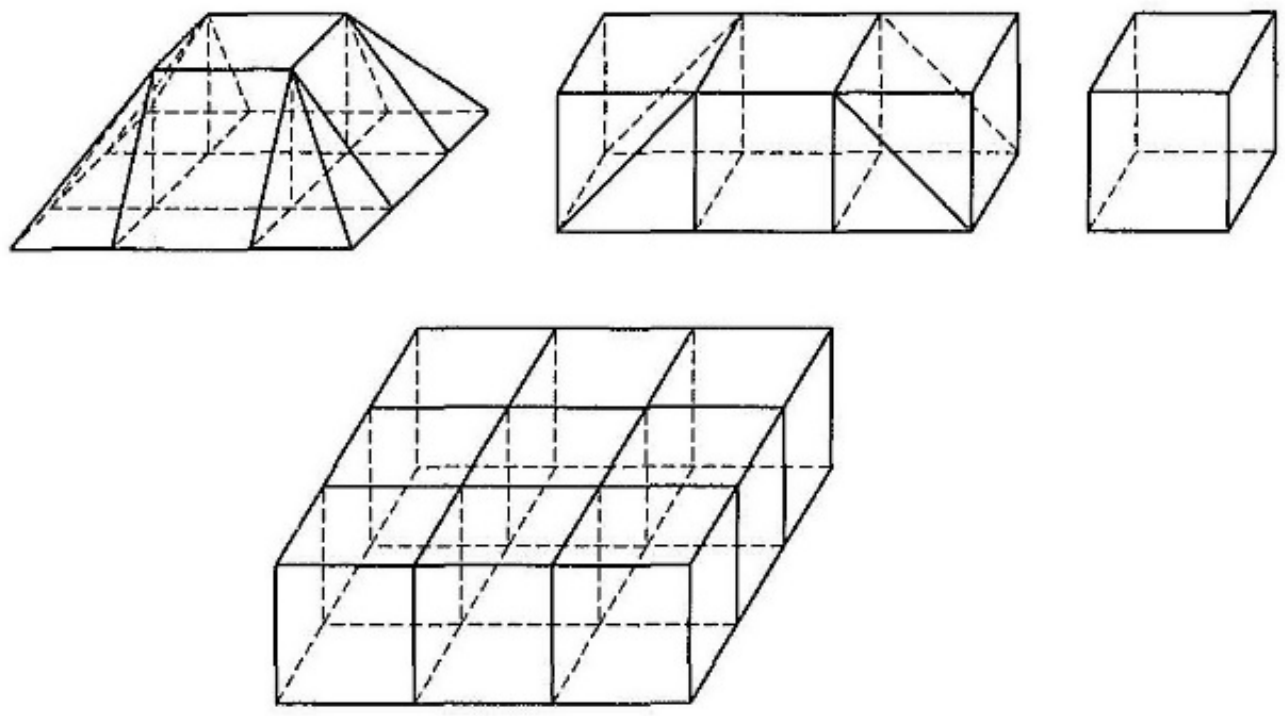

**Figure 5.13** – Traitement de la pyramide tronquée à base carrée par les blocs.

**Figure 6:** Geometrical representation of formula $F_T$ for the specific example of Liu Hui, with the volume of the frustum being one third of the sum of the three boxes which have to be filled following the classical formula. The figure taken from (Chemla/Guo 2004, 817).

Decisive for the understanding of his procedure is that Liu Hui did not choose an *arbitrary* way to combine the standard blocks but stipulated as the second and final step of his algorithm a very definite one:

The five cubes in the middle of the big box of nine cubes (marked in the figure below by diagonals on their faces) should be composed from one central cube and the 8 prisms (half cubes) of the two additional frustums, leaving the other remaining central cube as the third small box (Wagner 1979, 170). The four corner cubes in the big box would then result from integrating the total of 12 right asymmetrical pyramids (yangma) from the three frustums.

[36] Liu Hui always follows the classical formula of the Nine Chapters, he does not invent the boxes!

From our perspective it is immediately evident that in the general case with arbitrary $a, b, h$ these five cubes could be replaced by two crossing boxes of dimension $hab$ with a cuboid $hb$ in the middle thus leaving four corner cuboids with volume $h\left(\frac{a-b}{2}\right)^2$. Thus the ratio a:b = 3 : 1 is purely accidental and the big box with the quadratic base $a^2$ could be symmetrically stretched in accordance with a totally arbitrary $a$:$b$ along the marked area in figure 6a. But this "stretching" would keep intact the quadratic bases both in the central cuboid $h\,b^2$ and in the four corner cuboids $h\left(\frac{a-b}{2}\right)^2$.

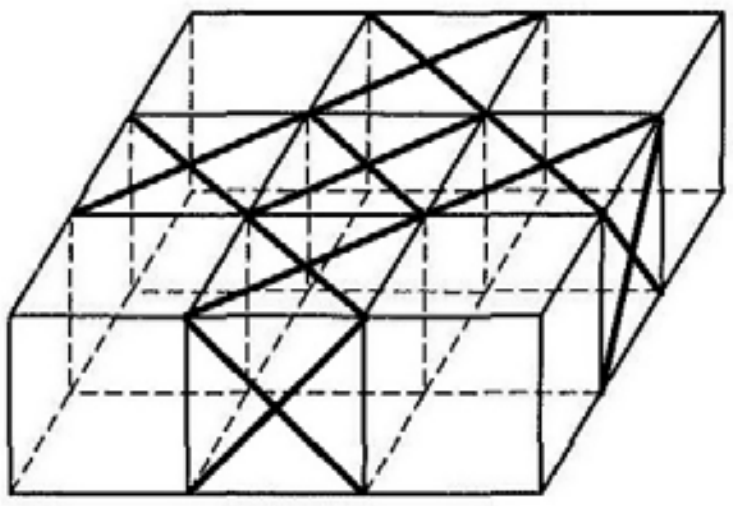

**Figure 6a:** The five cubes of the big box of formula $F_T$ (marked by diagonals) which Liu Hui composes from one central cube and 8 prisms (half cubes) of the two additional frustums

Despite of this arbitrariness of dimensions sweeping judgment that Liu Hui's proof is "general" is not immediately convincing.[37] To accept his proof as a general mathematical-geometrical conclusion requires a kind of mathematical abstraction from geometrical figure (congruence) to volume, which is comparable to the generalization of "equality" of plane figures (for instance triangles) from congruence to equality of areas in book I, 35 and from congruence to equality of volumes in book XI of Euclid's *Elements*.

At other places of his commentary Liu Hui says explicitly that a cuboid $ha^2$ with a height $h$ different from $a$ cannot usually be dissected into three asymmetrical "pyramids." Alluding to this insight of Liu Hui's, Chemla and Guo – on the basis of a much more careful analysis than can be presented here – come to the conclusion:

[37] I believe Van der Waerden glosses too quickly over the lack of explicit explanation within Liu Hui's commentary, when he says: "He works out the proof in the special case $h = 1$, $a = 1$, $b = 3$ only, but the idea of the proof is perfectly general." (Waerden 1983, 42) Van der Waerden does not even mention that these special values imply that the corner pyramids are inscribed in cubes which is not the general case. Likewise, he does not mention that the boxes produced from the three frustums are all composed of elementary cubic boxes. Only a full presentation of the restrictions of the example enables a discussion as to whether it is nevertheless general.

> "We can therefore see that Liu Hui cannot be ascribed the statement that the standard block verification method would be effective in solving all polyhedron volume problems, any more than he can be credited with the plan to develop his volume theory on this basis."[38]

Three facts, in particular, convince us that Liu Hui's proof with standard blocks can nevertheless be considered to be general. However, one needs to consider the broader context of his commentary and not just his text explaining problem 10 in chapter 5.

*First, and as indicated above*: the special dimensions chosen by Liu Hui do not affect the algorithm at least for the five central boxes (indicated with diagonals in figure 6a) within $ha^2$ which can be arbitrarily stretched to cuboids $hb^2$ and $h\left(\frac{a-b}{2}\right)b$.

*Second*: As far as the special dimensions of the corner pyramids are concerned, Liu Hui did have a proof for the general validity of $\frac{1}{3}$ – formula $F_P$ elsewhere in his commentary, using some geometrical operations similar to infinite series (Wagner 1979, 173), which however, probably went beyond the mathematical possibilities of the Egyptians. Liu Hui's was a method which is certainly based on "infinitesimal considerations" in the very broad sense as used by Neugebauer (see above). The physical realization of Liu Hui's method by using arbitrary models (provided he tried something along these lines) might cause some problems depending on the material used for the models. But one might imagine using 4 empty cuboid containers with height $h$ at the corners of the large cuboid $ha^2$ and filling them with material from the 12 corner pyramids of the three frustums after demolishing the latter models.

*Third*: Liu Hui gives at the end of his commentary an alternative procedure for calculating the volume of the frustum, where he – crucially – no longer talks about special measurements (elementary cubes). This is the alternative procedure mentioned above, which does not lead Liu Hui to the formula $F_T$ from the "Nine Chapters" but instead to the formula $F_{TA}$ which we know from our discussion above of the Egyptian case. Here the need and the conviction to have the $\frac{1}{3}$ – formula $F_P$ is even more visible. Liu Hui describes it verbally in the following manner:

[38] "On constat donc qu'on ne peut attribuer à Liu Hui la these selon laquelle la méthode de verification par blocs standard serait efficace pour résoudre tous les problems de volume de polyèdres, par plus qu'on ne peut lui prêter le projet de volouir, sur cette base, établir sa théorie des volumes." (Chemla/Guo 2004, 394)

"Another method: multiply by itself the difference of [the sides of] the squares, multiply by the height, and divide by 3. This gives [the volume of] the four yang-ma [corner pyramids]. Then multiply together [the sides of] the upper and lower squares and multiply by the height. This gives [the volume of] the central cube and the chien-tu [prisms] at the four sides. Add [these two results] to obtain the volume of the fang-ting [frustum]." (Wagner 1979, 170)

Here Liu Hui says nothing else than that the four asymmetrical corner pyramids with base edge $(a\text{-}b)/2$ of the frustum can be combined to one symmetrical pyramid with base edge $(a\text{-}b)$. Together with the parallelepiped box $hab$, which Liu Hui had found already quite at the beginning of his commentary (see above), this results in formula $F_{TA}$ which also Neugebauer, Gunn and Peet would derive around 1930, the only difference being that Liu Hui does not take out the common factor $h$.[39] Unlike his algorithm to confirm the classical formula Liu Hui does not connect his explanation to a calculation with standard blocks. This would of course result in the same volume, namely $\frac{1}{3} \bullet 1 \bullet (3-1)\,(3-1) + 1 \bullet 3 \bullet 1 = \frac{4}{3} + 3 = [= \frac{13}{3}]$. But the same calculation can be done with arbitrary $h$, $a$, and $b$. Liu Hui seems to say: we do not need the detour through the large cuboid $ha^2$ which leads to $F_T$ using three copies of the frustum on our way. We can alternatively calculate the volume more directly, using *only one* copy of the frustum, ignoring the revered formula in the "Nine Chapters." Except for also producing the parallelepiped box $hab$, one of the three boxes in the classical formula, the alternative algorithm is independent of Liu Hui's derivation of $F_T$.[40]

But this does not mean that Liu Hui disparages his own previous derivation of the classical formula $F_T$. His derivation in its general meaning (with general $a$, $b$ and $h$ as argued above) and not necessarily assuming the use of Liu Hui's "blocks" is visualized by the following picture taken from (Shutler 2009). The arrows (one here corrected) indicate in particular, that all 12 corner pyramids of the three copies of the truncated pyramid are supposed to go into the big box $ha^2$ as assumed in Liu Hui as well.

---

[39] The following remark is revealing here: "It seems that transparency of statement has been preferred before efficiency of computation." (Chemla 1990, 87).

[40] In (Chemla 1991, 84-87) the alternative procedure of Liu Hui is interpreted as a potential basis for another proof of the classical formula $F_T$ while I argue that Liu Hui emphasizes here his independence from the Nine Chapters. While I have problems to understand her interpretation at this point, I hasten to add that professional historians of Chinese mathematics such as Karine Chemla have added more reasons for the acceptance of the generality of Liu Hui's proofs, based on linguistic analysis. This goes beyond the more popular aims of this article.

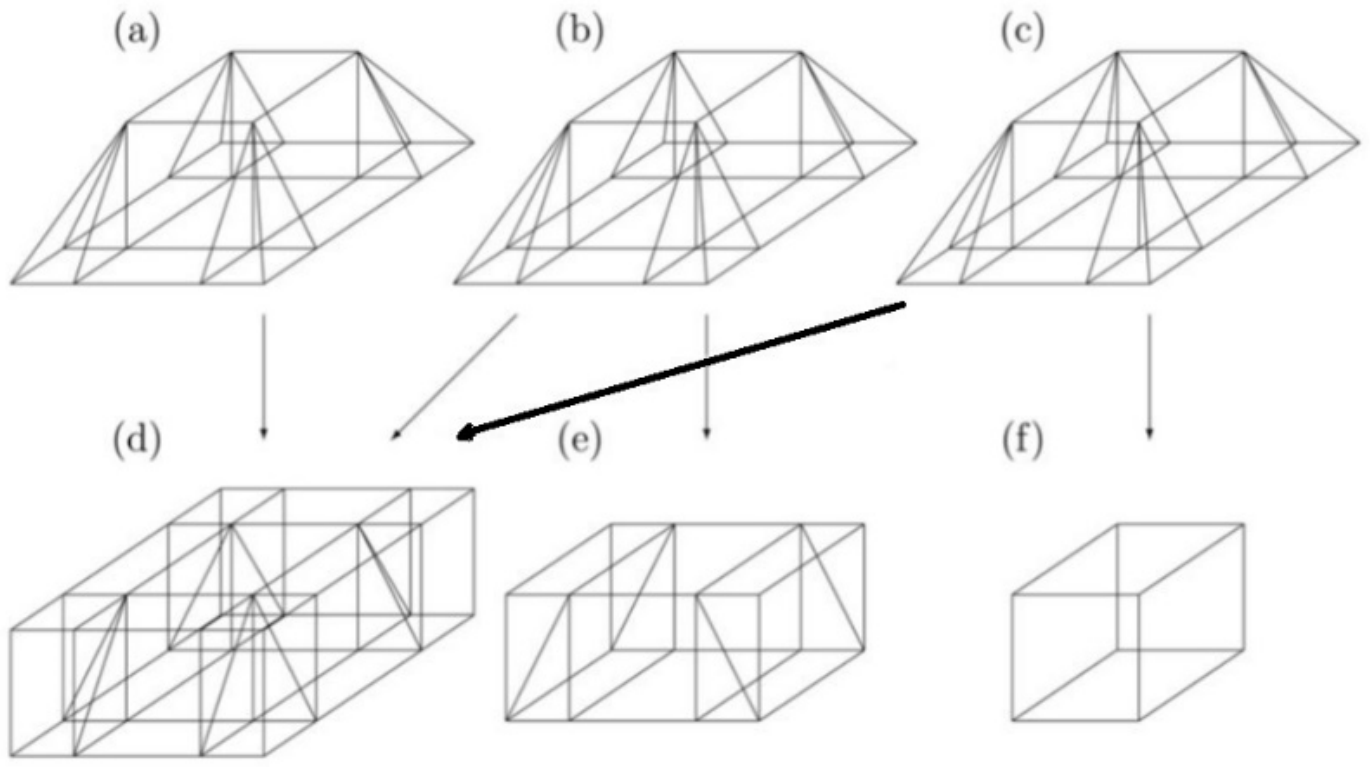

**Figure 2 (repeated from above):** from Shutler (2009, 349).

Note that unlike figure 6 and 6a the central 5 boxes within box (d) have arbitrary dimensions and the box in the middle is a cuboid $hb^2$. Shutler connects his figure immediately to the discussion of a possible Egyptian proof and does not mention the Chinese in his paper at all. I would assume though that – similar to my own experience – Shutler had Liu Hui's proof in the back of his mind after at least cursory encounters with literature on Chinese mathematics.

7. **An afterthought: How could the Egyptians (and the Chinese) have found the formula for the truncated pyramid in the first place?**

After summing up Liu Hui's proof, what does this tell us about a hypothetical proof of $F_T$ within Egyptian mathematics? And how could the Egyptians and the Chinese have found the formula in the first place, the latter being a heuristic question which differs from just looking for a possible proof (see section 2).

In fact, it comes as a kind of surprise that the three copies of the truncated pyramid produce exactly the three simple boxes $ha^2$, $hab$ and $hb^2$, whose volume can be easily calculated. One has probably to take into account the highly practical context both in Egypt and in China. Playing with models might have been encouraged or even been triggered by the "mathematical idea" that dissecting exactly *three* copies of the frustum would in any case *guarantee* the production of simple volumes because the 12 corner pyramids would then equal 4 cuboids in volume. But the final result, the small number of only three simple boxes, was in a way "coincidental". The resulting "symmetry" and relative "simplicity" of the formula $F_T$ was in my opinion not predictable. Admittedly, this is reminiscent of the typical unpredictability of formulas resulting from manipulations in algebra, and in this very restricted sense the Egyptians may already have possessed a certain "algebraic thinking," although performed with models (mentally or physically) and certainly without symbols. This is the core of what I consider the most convincing hypothetical way how the formula $F_T$ both in the Moscow Papyrus and in the Chinese Nine Chapters may have been originally found.

It remains unaccounted for why both the Egyptians and the Chinese found the "detour" through the large box $ha^2$ attractive and were not satisfied by the simple formula $F_{TA}$ which could be most easily derived from just one frustum. Regardless how the Egyptians (and later the Chinese) convinced themselves of the validity of the formula $F_T$, the existence of the algorithm of the Moscow papyrus seems to reveal a certain "theoretical" quality of Egyptian mathematics insofar as the rather elaborate algorithm is preferred to the geometrically more immediate, though less symmetrical formula $F_{TA}$.

If the "Egyptian or Chinese way of discovery" of the formula $F_T$ was such as conjectured above – namely tentatively looking at three copies of the frustum – a proof in the spirit of Liu Hui was almost a replication of that way of discovery. However, one should add "almost", because such proof was obviously based on the a priori information coming from the formula. In any case it required the conviction of the validity of the $\frac{1}{3}$ – formula $F_P$ if one should accept Liu Hui's explanation as a "proof" of $F_T$. Much depends on the "degree" of that conviction. In the Egyptian case I have argued that the conviction may have come from experience/experiment, in the Chinese case there have been, in addition, mathematical arguments (infinite series in geometric form). One should, however, be rather sure that Liu Hui's sophisticated "infinitesimal considerations" to prove $F_P$ (which have not been discussed in this paper) would have been beyond the methods available to the Egyptians.

Shutler, however, gives us a convincing argument, that an intuitive understanding of "infinitesimal considerations" of a more special kind might have been in reach of the Egyptians as well and may have contributed mathematically to their conviction of the $\frac{1}{3}$ – formula $F_P$.

## 8. Shutler's intuitive proof (2009) of the volume formula for the full pyramid

As mentioned in the beginning we have not even documentation that the Egyptians knew how to dissect a cube into three right asymmetrical pyramids or into six Juel pyramids. But how could they have even guessed the validity of the $\frac{1}{3}$ – formula $F_P$ for more general types of pyramids?

It was a pleasant surprise to me when my friend and colleague June Barrow-Green (London) drew my attention to Paul Shutler's 2009 article entitled "The problem of the pyramid or Egyptian mathematics from a postmodern perspective". First I looked at it with prejudice, because I had recently read too much about postmodernism, science and mathematics, and this without real benefit. I fear that more readers may have been discouraged by this packaging from reading Shutler's valuable paper.

Slowly it became clear to me that by "postmodern" the author seemed to mean in reality "premodern perspective" and that his work – without having historical research as a goal for itself – was written against "presentist" misinterpretations of history. As a mathematics didactician, Shutler was obviously primarily interested in using history and historical hypotheses to the benefit of mathematics teaching.

I read through the paper and had no judgment on Shutler's remarks on the Egyptian Rhind Papyrus. But then, in his remarks about the Moscow Papyrus and the truncated pyramid, I came across a simple and brilliant idea of Shutler's which aroused some envy in me.

He realized that in the special numerical example of the truncated pyramid in the Moscow Papyrus, where $a = 2b$, the missing pyramid at the top is congruent with the pyramid composed of the four corner pyramids, provided the figure is symmetrical. Because the entire pyramid is built on a square with twice the side length of the top pyramid one can conclude that its volume is eight times the volume of the top pyramid. (I discuss below to what extent one can assume that the Egyptians were able to understand this fact which follows from similarity and some simple infinitesimal considerations).[41] We have now two different relations between the volumes of the two different pyramids, the smaller and the big one, the known volume of the parallelepiped box *hab* mediating between the two. These two different relations allow a simple calculation of $F_P$.

In some more detail: two times the top pyramid $V_P$ plus the remaining parallelepiped box *hab* equal in volume eight times the top pyramid $V_P$. Thus the box has the same volume as six times the top pyramid which has been chopped off. With $a = 2b$ follows $6\ V_P = hab = 2hb^2$. This means that the volume of a pyramid $V_P$ with base $b^2$ is always equal to $\frac{h}{3}\ b^2$ because $b$ can be arbitrarily chosen.

[41] (Thomas 1931, 50) mentions the 1/8 ratio too and assumes the Egyptians could have had this insight as well. But he does not draw a comparison to the corner pyramids in the truncated pyramid, while Shutler does.

Again, this could be interpreted in a "presentist" way as non-symbolic "algebra", namely as "two independent linear equations". A more intuitive understanding, however, would be that this is a simple comparison of geometrical figures, for example by manipulation with models of parts of the pyramid, something which is not easily comparable to the modern manipulation of algebraic equations and symbols. It should therefore have been within the reach of the Egyptians. Shutler assumes that the concrete (rather moderate) dimensions given in the papyrus with $a = 2b$ indicate that the Egyptians at least occasionally used exactly such models, which would then have shown the congruence of the missing top pyramid with the combination of the smaller asymmetrical corner pyramids. This does not seem to contradict the assumptions of Egyptologists either.

Shutler calls this a "proof by example", by which he means that the Egyptians were able "to derive the general case from a certain numerical example" (Shutler 2009, 350), and recommends such a proof for mathematics teaching. This would imply a much broader interpretation of problem 14 in the Moscow Papyrus than hitherto accepted. Until now, the text has mostly been understood as a general "claim by example", and it has been interpreted as a general formula without proof.

I hesitate to follow Shutler on this point, above all because the numerical example must be supplemented by a second, quite different insight, if one wants to derive a proof from it. This is the ratio of the volumes of the top pyramid to the complete pyramid of 1:8, which is obviously based on what Neugebauer would have called "infinitesimal considerations".

But how can we assume that the Egyptians had this "second insight"? Shutler only mentions time as a factor which he thinks would have shown that for half of the height of the full pyramids seven eighths of the total construction time were required, an argument which I find less convincing than a similar one related to the material needed.[42] But of course, one could imagine that the Egyptians found the ratio by using physical models of both the top pyramid and of the full pyramid and weighing them or measuring their volumes within water.

[42] Exactly the time spent may not be a reliable marker because the upper part of the pyramid may take much more time to build due to the need to transport the bricks to the higher level.

But arguing more mathematically (geometrically): “Similarity” at least in an intuitive, non-formalized sense seems to be decisive, although Shutler does not explicitly say so. Undoubtedly the Egyptians knew, both practically and more “theoretically” by geometrical imagination, that a cube with twice the side of a smaller cube has eight times the volume of the smaller. The next logical step is to assume that two similar three-dimensional figures – by similarity meaning the intuitive notion of a constant stretching factor in all directions – can be filled by more and more diminishing cubes in the volume relation 1 : 8.

At first it looks as though one needs some deeper understanding of similarity in order to see that the height of the pyramid changes in the same proportion as the side-length of the top square. But for the special case $a = 2b$ we actually only need a very limited understanding of similarity, in particular no similarity proportions, as the following simple illustration, based on the congruence of two triangles, immediately shows:

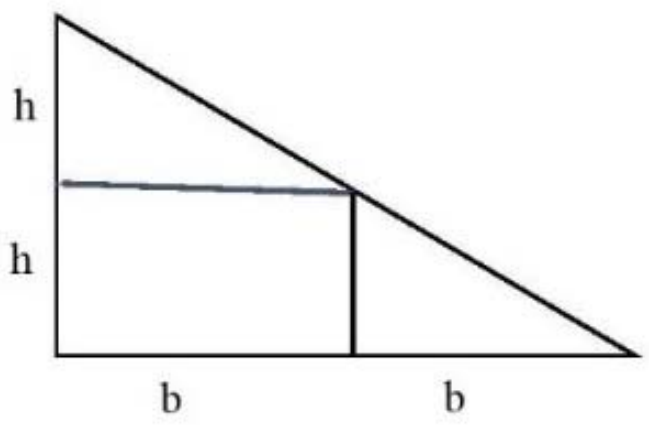

**Figure 7:** The picture shows trivially that $a = 2\,b$ (here for the asymmetrical pyramid) implies the doubling of the height $h$ for the full pyramid

As to a more sophisticated use of similarity by Egyptian mathematicians than assumed here one has to take into account that even the notion of a “triangle” in our modern sense cannot be uncritically assumed to have existed in ancient Egypt.[43]

## 9. Summary, conclusions, and open questions

Modern mathematical results often enable a fresh look at history and explain limitations and struggles of mathematicians of the past. This is the case for instance for the proof of the impossibility of squaring the circle and for the discovery of non-Euclidean geometries, both in the 19th century. Dehn's solution to Hilbert's Third Problem (1900) which was alluded to above allows conclusions as to which finite dissections of polyhedra are logically and mathematically possible and can thus serve the comparison of volumes.

[43] Comparable restrictions are even valid for Chinese mathematics 2000 years later: see (Chemla/Guo 2004, 662). I leave it to the specialists to decide whether certain Egyptian algorithms such as the ‘sqd-calculation’ are based on implicit (!) similarity proportions in triangles,

But there are dangers connected to such retrospective look into history too. In hindsight it is surprising that in the almost 100 years from the first interpretation of the Moscow Papyrus (1917) until Shutler (2009), no interpreter with the possible exception of Van der Waerden (1983) (see below) seems to have thought of an Egyptian proof of the kind of Liu Hui's, i.e. multiplying the formula $F_T$ by 3 and interpreting it as a volume equality of three identical copies of the frustum with three easily calculable boxes:

$$V_T + V_T + V_T = ha^2 + hab + hb^2$$

After the widely accessible publication in (Wagner 1979) of Liu Hui's commentary on Chinese pyramids had appeared one should have expected comparisons with a hypothetical Egyptian proof. While mathematically minded Egyptologists around 1930 could not know the Chinese case, the sinologist Wagner did not mention the Moscow Papyrus for comparison, because he did not approach the problem from the standpoint of the general history of mathematics.[44] It is, however, somewhat surprising that historians of Egyptian mathematics after Wagner tended not to draw on Chinese mathematics for comparison.[45] One may surmise that they took unsuccessful speculations of their forbears around 1930 as methodological warning signs.

As a result of this, the non-historian Shutler (2009), who was mainly interested in Egyptian mathematics for pedagogical purposes, could not easily find Liu Hui or Wagner mentioned in the literature related to the Egyptian problem.[46] He probably even did not know about them or had Liu Hui only vaguely in the back of his mind. However, coming from mathematical didactics had an advantage too. Shutler gave a plausible argument that the $\frac{1}{3}$ – formula for the volume of the full pyramid was within the reach of Egyptian mathematicians and he "reinvented" afterwards Liu Hui's argument as a kind of matter of course (figure 2).

[44] Wagner confirmed in an email to the author, 25 March 2020, that he did not know the Egyptian formula at the time.
[45] Chinese mathematics is not mentioned in Imhausen (2016). More generally the author avoids references to other mathematical cultures, being aware of lack of actual historical interactions.
[46] Shutler mentions neither Chinese mathematics nor (Waerden 1983) where he could have found a stimulus for comparison.

Shutler's proposal could be called "exact intuition" as opposed to the "inexact" one which is conveyed by the classical formula $F_T$, as explained above. In any case this proof method can be easily made rigorous by using fundamental notions of modern mathematics.

A standpoint of his own had the prominent mathematician Van der Waerden (1983), Neugebauer's student in history of mathematics in 1920s' Göttingen. He was interested in all ancient mathematical cultures and drew ambitious, intercultural conclusions on an insufficient documentary basis even going as far as hypothesizing communication channels at the time. This approach has its pedagogical merits and parts of my paper are inspired by it because he was not afraid to draw parallels between disparate cultures such as the Egyptian and the Chinese. With respect to our problem Van der Waerden has a general conjecture about relations between Egyptian and Chinese mathematics, which goes too far by general agreement of historians:

> "One and the same correct rule for the volume of a truncated pyramid.
> In my opinion, this network of interrelations and similarities can only
> be explained by assuming a common origin for the mathematics and astronomy
> of these ancient countries." (Waerden 1983, 44)

It seems that strong and insufficiently founded conjectures like this discouraged other historians to venture the possibility that the Egyptians thought geometrically in a similar manner like the Chinese 2000 years later. The overall criticism which met Van der Waerden's book may have acted for professional historians as another call for increased caution.

Meanwhile there exists a new generation of historians of Babylonian, Greek, Chinese, Indian and Arabic mathematics, who can draw on a much broader documentary basis for these cultures than for the Egyptian, and who often turn to intercultural comparison (Chemla 2012).

Using the work of professional historians of Egyptian and Chinese mathematics, I have tried to argue that both Liu Hui's actual proof (ca 263 CE) and Shutler's hypothetical proof (2009) of the Egyptian/Chinese formulae for the volume of the truncated and for the volume of the full pyramid were within the reach of Egyptian mathematicians as well.

While Liu Hui's proof is intuitive and rigorous even in the modern sense, Shutler's proof for the volume of the full pyramid is intuitive and convincing though not rigorous but can be easily turned into a rigorous one by using elementary notions from infinitesimal analysis. Both are helpful by connecting to the intuitions of students today.

I have further argued in this paper that within the reach of Egyptian mathematicians there were probably other, alternative formulas for the volume of the truncated pyramid which were no more complicated than, and as easy to memorize as the one in the Moscow Papyrus which I called $F_T$. In my opinion, the formula in the Moscow Papyrus can even be considered to be more "abstract" than for instance a particular alternative formula $F_{TA}$ found by Liu Hui 2000 years later. Historians around 1930 have implicitly recognized the "abstract" character of $F_T$ when they emphasized the "symmetry" of that formula for the volume of the truncated pyramid. Some of them, however, went too far in the direction of a "presentist" methodological approach, and speculated about some geometric and even hidden algebraic thinking related to plane geometry among the ancient Egyptians, ignoring the possibility of more intuitive ways based on considering volumes directly. To put it a little more pointedly: The Egyptians were probably "better" mathematicians than can be documented with historical certainty, but not as "good" (modern) as some historians have assumed by backward projection of their own posterior knowledge.

Even if this interpretation goes too far, one may still assume that the Egyptians had already something of the playful mind of the modern mathematician who develops general formulas and enjoys their beauty and symmetry without always paying attention to their usefulness for applications.

While the real motives and methods of the Egyptians will probably – not least due to a lack of relevant documents – remain hidden under the veil of history, their formulas for the volume of the truncated and the full pyramid remain pieces of important elementary mathematics. A discussion of the hypothetical derivation of the two formulas, in particular the role of special cases and to which extent they point to general meaning, and the importance of experimental or of purely mathematical reasoning can be stimulating for today's teaching of mathematics.

**Acknowledgments**

I am grateful to the following colleagues and friends for reading earlier versions of this manuscript and giving helpful advice: June Barrow-Green (London), Reinhard Bölling (Berlin), Karine Chemla (Paris), Joseph Dauben (New York), Simon Goodchild (Kristiansand), Christopher Hollings (Oxford), Jesper Lützen (Copenhagen), Karl-Heinz Schlote (Altenburg), Nathan Sidoli (Tokyo), and Donald Wagner (Copenhagen).